**Preliminary results about a wide class of nonlinear recursive sequences**

M. De la Sen

Instituto de Investigacion y Desarrollo de Procesos. Universidad del Pais Vasco

Campus of Leioa (Bizkaia) – Aptdo. 644- Bilbao, 48080- Bilbao. SPAIN

**Abstract**. This paper deals with the boundedness of the solutions and the stability properties of the positive solutions of recursive sequence $x_{k+1} = A_k + \sum_{i=1}^{q} B_{ik} \frac{x_{k+1-\ell_i}^{p_i}}{x_{k+1-s_i}^{r_i}}$ where $A_k$ and $B_{ik}$ are positive real sequences, $p_i$ and $r_i$ are real constants and $\ell_i$ and $s_i$ are positive integers. The above recursive sequence generalizes some previously studied recursive sequences in the literature.

## 1. Introduction

The positive solutions of recursive sequence:

$$x_{k+1} = A_k + \sum_{i=1}^{q} B_{ik} \frac{x_{k+1-\ell_i}^{p_i}}{x_{k+1-s_i}^{r_i}} \ ; \ \forall k \in \mathbf{Z}_{0+} = \mathbf{Z}_+ \cup \{0\} \tag{1.1}$$

where $A_k$ and $B_{ik}$ are positive real sequences, $p_i$ and $r_i$ are real constants and $\ell_i$ and $s_i$ are positive integers; $\forall i \in \overline{q} := \{1, 2, \dots, q\}$. The above recursive sequence generalizes some previously proposed ones investigated in the background literature. It is investigated from a equilibrium stability and solution boundedness point of view. The motivation of the choice of the recursive sequence (1.1) is the significant interest which has been devoted in the recent past to recursive sequences of close structures where quotients of different delayed values of the solution have been introduced to construct their solutions. Some of those previous studies which have inspired this one are listed below:

1) Yang, Sun and Zhang investigate in [1] the equilibrium points, boundedness of the solutions and oscillatory character of the recursive difference equation:

$$x_n = A + \frac{x_{n-k}^p}{x_{n-1}^p} \ ; \ \forall n \in \mathbf{Z}_+$$

2) Amleh, Grove and Ladas study in [2] the second-order rational recursive difference equation:

$$x_{n+1} = \alpha + \frac{x_{n-1}}{x_n} \ ; \ \forall n \in \mathbf{Z}_{0+}$$

3) Berenhaut and Stevic study in [3] and [4] the subsequent more general difference equations:

$$x_{n+1} = \alpha + \frac{x_{n-k}^p}{x_n^p} \ ; \qquad x_n = A + \left( \frac{x_{n-2}}{x_{n-1}} \right)^p$$

4) Stevic investigates the difference equations below:

$$x_{n+1} = \alpha + \frac{x_{n-1}^p}{x_n^p} \ ; \ x_{n+1} = \alpha + \frac{x_n^p}{x_{n-1}^p} \ ; \ x_{n+1} = \alpha + \frac{x_n^p}{x_{n-1}^r} : \ \forall n \in \mathbf{Z}_+$$

in [5], [6] and [7], respectively.

5) Dehghan and Mazrooei-Sebdani investigate in [8] the properties of stability, periodicity, boundedness and oscillatory behaviour of the rational difference equation:

$$x_{n+1} = \frac{\alpha + \gamma x_{n-1}}{A + Bx_n}$$

Furthermore, it has to be pointed out that rational particular form of the difference equation in [1] for the parameter r being unity, which is also close to particular forms of the rational difference equation proposed in [3], has been investigated in [14]. The stability properties and boundedness of the solutions of a class of rational difference equations with periodic coefficients and those of second order systems of difference equations have been investigated in [13] and [15], respectively. On the other hand, the global stability and asymptotic of continuous- time dynamic systems under point and distributed delayed dynamics has been investigated in [16] and references therein. The provided results can be a starting point of parallel results concerning the discretization of such dynamic systems.

**2. Boundedness of the solutions**

This paper first considers the solutions of the following difference equation:

$$x_{k+1} = A_k + B_k \frac{x_{k+1-\ell}^{p}}{x_{k+1-s}^{r}}\ ; \ \forall k \in \mathbf{Z}_{0+} \tag{2.1}$$

which is a particular case of the recursive sequence (1.1) under admissible initial conditions:

$$(x_{-j} > 0;\ \forall j \in [0, s] \cap \mathbf{Z} \text{ if } s \ge \ell)$$

$$\wedge [(x_{-j} > 0;\ \forall j \in [0, s] \cap \mathbf{Z}) \wedge (x_{-j} \ge 0;\ \forall j \in [s+1, \ell] \cap \mathbf{Z}) \text{ if } s < \ell]$$

with given constants $p \in \mathbf{R}$, $r \in \mathbf{R}$, $\ell \in \mathbf{Z}_+$, $s \in \mathbf{Z}_+$ and real positive uniformly bounded sequences $\{A_k\}_{k \in \mathbf{Z}_{0+}}$ and $\{B_k\}_{k \in \mathbf{Z}_{0+}}$, respectively.

The following result of uniform boundedness of the solution sequence of the discrete equation (2.1) follows:

**Theorem 2.1**. The following properties hold:

**(i)** Assume that there exist constants m, $M(\ge m)$, $a_m$, $a_M(\ge a_m)$, $b_m$, $b_M(\ge b_m) \in \mathbf{R}_+$ which satisfy the subsequent constraints:

$$a_m \le \min_{0 \le k < \infty} A_k \le \max_{0 \le k < \infty} A_k \le a_M\ ; \quad b_m \le \min_{0 \le k < \infty} B_k \le \max_{0 \le k < \infty} B_k \le b_M \tag{2.2}$$

$$a_m \le m \le \inf_{\min(1-\ell, 1-s) \le j \le 0} x_j \le \sup_{\min(1-\ell, 1-s) \le j \le 0} x_j \le M \tag{2.3}$$

$$m\left(1 - b_m \frac{m^{p-1}}{M^r}\right) \le a_m \le a_M \le M\left(1 - b_M \frac{M^{p-1}}{m^r}\right) \tag{2.4}$$

$$M \geq \max\left( m, b_m^{1/r} m^{\frac{p-1}{r}} \right); \quad m \geq b_M^{1/r} M^{\frac{p-1}{r}} \tag{2.5}$$

$$\left( \frac{b_M}{b_m} \right)^{p-1} \geq \frac{M}{m} \geq 1 \tag{2.6}$$

Then, $m \leq \inf_{\min(1-\ell, 1-s) \leq k < \infty} x_k \leq \sup_{\min(1-\ell, 1-s) \leq k < \infty} x_k \leq M$. (2.7)

**(ii)** Property (i) can be extended directly for the case $b_m \in \mathbf{R}_{0+}$, $b_M \in \mathbf{R}_+$, $\min_{0 \leq k < \infty} B_k < \max_{0 \leq k < \infty} B_k$ by removing the constraint (2.7) and modifying the second inequality in the second relation (2.2) to a strict one.

**Proof**: Proceed by complete induction by assuming that

$$m \leq \inf_{\min(1-\ell, 1-s) \leq j \leq i} x_j \leq \sup_{\min(1-\ell, 1-s) \leq j \leq i} x_j \leq M; \ \forall i \in \overline{k}, \text{ some given } k \in \mathbf{Z}_+ \tag{2.8}$$

Assume that $\min(p, r) \geq 0$. Then, one gets from (2.1) that

$$\left[ m \leq \inf_{\min(1-\ell, 1-s) \leq j \leq i} x_j \leq \sup_{\min(1-\ell, 1-s) \leq j \leq i} x_j \leq M \right]; \ \forall i \in \overline{k} := \{1, 2, \ldots, k\} \Rightarrow m \leq x_{k+1} \leq M \,]$$

$$\Leftrightarrow \left[ m \leq \inf_{\min(1-\ell, 1-s) \leq j \leq k+1} x_j \leq \sup_{\min(1-\ell, 1-s) \leq j \leq k+1} x_j \leq M \right] \tag{2.9}$$

Assume that the constraint (2.4) holds , subject to the necessary condition (2.5) (which in turn requires the necessary condition (2.6)) in view of (2.1) subject to (2.2)-(2.3). Thus, (2.9) guarantees (2.4) via complete induction. The proof of Property (i) is complete. Property (ii) follows since Property (i) can be extended for the case $b_m = 0$ with the given modifications. □

Theorem 2.1 might be simplified by noting that $a_m \leq \inf_{\min(1-\ell, 1-s) \leq j \leq 0} x_j$ from (2.1). Thus, $a_m$ may be chosen replacing m as the uniform lower-bound of all the elements of the sequence $\{x_k\}$ under a sufficiently small $b_M$ related to an explicit upper-bound obtained from M and $a_m$. This philosophy is formally addressed in the subsequent result:

**Corollary 2.2**. Assume that there exist real constants M, $a_m, b_M \in \mathbf{R}_+$ which satisfy the constraints:

$$a_m \leq \min_{0 \leq k < \infty} A_k \leq \max_{0 \leq k < \infty} A_k \leq M \left( 1 - b_M \frac{M^{p-1}}{a_m^r} \right)$$

$$0 < \min_{0 \leq k < \infty} B_k \leq \max_{0 \leq k < \infty} B_k \leq b_M \leq \frac{a_m^r}{M^{p-1}}$$

$$a_m \le \inf_{\min(1-\ell,1-s)\le j\le 0} x_j \le \sup_{\min(1-\ell,1-s)\le j\le 0} x_j \le M \tag{2.10}$$

Then, (2.7) holds with the replacement $m \to a_m$. The result also holds for the case $\min_{0\le k<\infty} B_k \ge 0$, $\max_{0\le k<\infty} B_k > 0$.

**Proof**: Proceed by complete induction by assuming that (2.8) holds with then replacement $m \to a_m$. Then, one gets from (2.1), subject to the constraints (2.10), that (2.9) also holds with then replacement $m \to a_m$. As a result, (2.7) also holds with such a replacement . □

The set of constraints required by Corollary 2.2 is less involved than those of Theorem 2.1. However, the price to be paid is that the uniform lower-bound estimate to be used for the solution sequence of (2.1) is necessarily $a_m$ and never a potential larger real lower-bound. Theorem 2.1 and Corollary 2.2 have the following simpler parallel result. All the elements of the sequence $\{x_k\}$ belong to the interval $[m, M] \subset \mathbf{R}_+$, provided that the initial conditions of (2.1) belong to such an interval, if all the member of the parameterizing sequence $\{A_k\}$ belong to some prescribed interval in $\mathbf{R}_{0+}$ provided that $B_k \in [b_m, b_M] \subset \mathbf{R}_{0+}$ for some $b_m \in \mathbf{R}_{0+}$, $b_M (\ge b_m) \in \mathbf{R}_+$.

**Corollary 2.3**. Assume that there exist $m, M(\ge m) \in \mathbf{R}_+$, $b_m \in \mathbf{R}_{0+}$, $b_M (\ge b_m) \in \mathbf{R}_+$ which satisfy:

$$m\left(1 - b_m \frac{m^{p-1}}{M^r}\right) \le A_k \le M\left(1 - b_M \frac{M^{p-1}}{m^r}\right) \quad ; \quad B_k \in [b_m, b_M] \subset \mathbf{R}_{0+} \tag{2.11}$$

subject to the necessary conditions (2.5)-(2.6). If the initial conditions are constrained to:

$$m \le \inf_{\min(1-\ell,1-s)\le j\le 0} x_j \le \sup_{\min(1-\ell,1-s)\le j\le 0} x_j \le M \tag{2.12}$$

Then, (2.7) holds.

**Proof**: It follows by complete induction in a similar way as the proof of Theorem 2.1. □

The sequence (2.1) is generalized to the subsequent one whose study is the main interest of this manuscript:

$$x_{k+1} = A_k + \sum_{i=1}^{q} B_{ik} \frac{x_{k+1-\ell_i}^{p_i}}{x_{k+1-s_i}^{r_i}} \; ; \; \forall k \in \mathbf{Z}_{0+} \tag{2.13}$$

under admissible initial conditions:

$$(x_{-j} > 0;\ \forall j \in [0, s] \cap \mathbf{Z} \text{ if } s \ge \ell)$$

$$\wedge [(x_{-j} > 0;\ \forall j \in [0, s] \cap \mathbf{Z}) \wedge (x_{-j} \ge 0;\ \forall j \in [s+1, \ell] \cap \mathbf{Z}) \text{ if } s < \ell]$$

with given constants $p \in \mathbf{R}$, $r \in \mathbf{R}$, $\ell \left( := \max_{i \in \bar{q}} \ell_i \right) \in \mathbf{Z}_+$, $s \left( := \max_{i \in \bar{q}} s_i \right) \in \mathbf{Z}_+$ and real positive uniformly bounded sequences $\{A_k\}_{k \in \mathbf{Z}_{0+}}$, $\{B_{ik}\}_{k \in \mathbf{Z}_{0+}}$; $\forall i \in \bar{q}$, respectively. For the estimation of uniform upper and lower bounds of the solution sequence, define the following real constants:

$$\mu_p := \left( p_j : M^{p_j - 1} := \max\left( M^{p_i - 1} : i \in \bar{q} \right), \text{ some } j \in \bar{q} \right)$$

$$\nu_r := \left( r_j : M^{r_j} := \min\left( M^{r_i} : i \in \bar{q} \right), \text{ some } j \in \bar{q} \right)$$

$$\nu_p := \left( p_j : m^{p_j - 1} := \min\left( m^{p_i - 1} : i \in \bar{q} \right), \text{ some } j \in \bar{q} \right)$$

$$\mu_r := \left( r_j : m^{r_j} := \max\left( m^{r_i} : i \in \bar{q} \right), \text{ some } j \in \bar{q} \right) \tag{2.14}$$

for some real constants m, $M(\geq m) \in \mathbf{R}_+$ such that:

$$a_m \leq m \leq \inf_{\min(1-\ell, 1-s) \leq j \leq 0} x_j \leq \sup_{\min(1-\ell, 1-s) \leq j \leq 0} x_j \leq M \tag{2.15}$$

The particular coefficients $p_i$, $r_i$ $(i \in \bar{q})$ which define the definitions (2.14) depend on the coefficients $p_i$ and $r_i$ being non less than unity and non less than zero, respectively , if any and also if the real constants M and m are less than unity or not. Related to these considerations, the subsequent technical result follows by direct inspection of (2.14) .

**Lemma 2.4**. Decompose the sets of the coefficients $p_i, r_i$ be $\Lambda_p := \{p_i : i \in \bar{q}\}$ and $\Lambda_r := \{r_i : i \in \bar{q}\}$ uniquely as the respective disjoint unions of the sets $\Lambda_{p+} := \{p_i \geq 1 : i \in \bar{q}\}$ and $\Lambda_{p-} := \{p_i < 1 : i \in \bar{q}\}$; and $\Lambda_{r+} := \{r_i \geq 0 : i \in \bar{q}\}$ and $\Lambda_r := \{r_i < 0 : i \in \bar{q}\}$.The following properties hold:

**(i.1)** $\left( \Lambda_{p+} \neq \varnothing \wedge M \geq 1 \right) \Rightarrow \mu_p := \left( \max\ p_j : p_j \in \Lambda_{p+} \right)$
**(i.2)** $\left( \Lambda_{p+} \neq \varnothing \wedge \Lambda_{p-} = \varnothing \wedge M < 1 \right) \Rightarrow \mu_p := \left( \min\ p_j : p_j \in \Lambda_{p+} \right)$
**(i.3)** $\left( \Lambda_{p+} \neq \varnothing \wedge \Lambda_{p-} \neq \varnothing \wedge M < 1 \right) \Rightarrow \mu_p := \left( \min\ p_j = -\max \left| p_j \right| : p_j \in \Lambda_{p-} \right)$
**(i.4)** $\left( \Lambda_{r+} \neq \varnothing \wedge M \geq 1 \right) \Rightarrow \nu_r := \left( \min\ r_j : r_j \in \Lambda_{r+} \right)$
**(i.5)** $\left( \Lambda_{r+} \neq \varnothing \wedge \Lambda_{r-} = \varnothing \wedge M < 1 \right) \Rightarrow \nu_r := \left( \max\ r_j : r_j \in \Lambda_{r+} \right)$
**(i.6)** $\left( \Lambda_{r+} \neq \varnothing \wedge \Lambda_{r-} \neq \varnothing \wedge M < 1 \right) \Rightarrow \nu_r := \left( \max\ r_j = -\min \left| r_j \right| : r_j \in \Lambda_{r-} \right)$

**(ii.1)** $\left( \Lambda_{p+} \neq \varnothing \wedge m \geq 1 \right) \Rightarrow \nu_p := \left( \min\ p_j : p_j \in \Lambda_{p+} \right)$
**(ii.2)** $\left( \Lambda_{p+} \neq \varnothing \wedge \Lambda_{p-} = \varnothing \wedge m < 1 \right) \Rightarrow \nu_p := \left( \max\ p_j : p_j \in \Lambda_{p+} \right)$
**(ii.3)** $\left( \Lambda_{p+} \neq \varnothing \wedge \Lambda_{p-} \neq \varnothing \wedge m < 1 \right) \Rightarrow \nu_p := \left( \max\ p_j = -\min \left| p_j \right| : p_j \in \Lambda_{p-} \right)$
**(ii.4)** $\left( \Lambda_{r+} \neq \varnothing \wedge m \geq 1 \right) \Rightarrow \mu_r := \left( \max\ r_j : r_j \in \Lambda_{r+} \right)$
**(ii.5)** $\left( \Lambda_{r+} \neq \varnothing \wedge \Lambda_{r-} = \varnothing \wedge m < 1 \right) \Rightarrow \mu_r := \left( \min\ r_j : r_j \in \Lambda_{r+} \right)$
**(ii.6)** $\left( \Lambda_{r+} \neq \varnothing \wedge \Lambda_{r-} \neq \varnothing \wedge m < 1 \right) \Rightarrow \mu_r := \left( \min\ r_j = -\max \left| r_j \right| : r_j \in \Lambda_{r-} \right)$

**(iii.1)** $\left(\Lambda_{p+} = \varnothing \wedge M \geq 1\right) \Rightarrow \mu_p := \left(\max p_j = -\min\left|p_j\right| : p_j \in \Lambda_{p-}\right)$

**(iii.2)** $\left(\Lambda_{p+} = \varnothing \wedge M < 1\right) \Rightarrow \mu_p := \left(\min p_j = -\min\left|p_j\right| : p_j \in \Lambda_{p-}\right)$

**(iii.3)** $\left(\Lambda_{r+} = \varnothing \wedge M \geq 1\right) \Rightarrow \nu_r := \left(\min r_j = -\min\left|r_j\right| : r_j \in \Lambda_{r-}\right)$

**(iii. 4)** $\left(\Lambda_{r+} = \varnothing \wedge M < 1\right) \Rightarrow \nu_r := \left(\max r_j = -\min\left|r_j\right| : r_j \in \Lambda_{r-}\right)$

**(iv.1)** $\left(\Lambda_{p+} = \varnothing \wedge m \geq 1\right) \Rightarrow \nu_p := \left(\min p_j = -\min\left|p_j\right| : p_j \in \Lambda_{p-}\right)$

**(iv.2)** $\left(\Lambda_{p+} = \varnothing \wedge m < 1\right) \Rightarrow \nu_p := \left(\max p_j = -\min\left|p_j\right| : p_j \in \Lambda_{p-}\right)$

**(iv.3)** $\left(\Lambda_{r+} = \varnothing \wedge m \geq 1\right) \Rightarrow \mu_r := \left(\max r_j = -\min\left|r_j\right| : r_j \in \Lambda_{r-}\right)$

**(iv.4)** $\left(\Lambda_{r+} = \varnothing \wedge m < 1\right) \Rightarrow \mu_r := \left(\min r_j = -\min\left|r_j\right| : r_j \in \Lambda_{r-}\right)$ □

The following result is the generalization of Theorem 2.1 if the discrete equation (2.1) is replaced by its generalization (2.13):

**Theorem 2.5**. Assume that the sequences $\left\{A_k\right\}$ and $\left\{B_{ik}\right\}$; $i \in \overline{q}$ satisfy:

$$0 < a_m \leq \min_{0 \leq k < \infty} A_k \leq \max_{0 \leq k < \infty} A_k \leq a_M$$

$$0 \leq b_{im} \leq \min_{0 \leq k < \infty} B_{ik} \leq \max_{0 \leq k < \infty} B_{ik} \leq b_{iM} \;;\; \max_{0 \leq k < \infty} B_{ik} > 0$$

$$m\left(1 - b_m \frac{m^{\nu_p - 1}}{M^{\mu_r}}\right) \leq a_m + b_m \frac{m^{\nu_p}}{M^{\mu_r}} \leq a_m \leq a_M \leq M\left(1 - b_M \frac{M^{\mu_p - 1}}{m^{\nu_r}}\right)$$

$$b_m := \sum_{i=1}^{q} b_{im} \;;\; b_M := \sum_{i=1}^{q} b_{iM}$$

$$m\left(1 - b_m \frac{m^{\nu_p - 1}}{M^{\mu_r}}\right) \leq a_m \leq a_M \leq M\left(1 - b_M \frac{M^{\mu_p - 1}}{m^{\nu_r}}\right) \tag{2.16}$$

$$M \geq \max\left(m, b_m^{1/\mu_r} m^{\frac{\nu_p - 1}{\mu_r}}\right) ;\; m \geq b_M^{1/\nu_r} M^{\frac{\mu_p - 1}{r}} \tag{2.17}$$

with the real constants (2.14) being calculated from Lemma 2.4, and that the solution sequence of (2.13) is subject to initial conditions satisfying:

$$m \leq \inf_{\min(1-\ell, 1-s) \leq j \leq 0} x_j \leq \sup_{\min(1-\ell, 1-s) \leq j \leq 0} x_j \leq M$$

Then, (2.7) holds.

**Proof**: It follows by complete induction. On gets directly from (2.13) and (2.16) that

$$m \leq a_m + b_m \frac{m^{\nu_p}}{M^{\mu_r}} \leq a_M + b_M \frac{M^{\mu_p}}{m^{\nu_r}} \leq M \tag{2.18}$$

with $\mu_p$, $\mu_r$, $\nu_p$ and $\nu_r$ defined in Lemma 2.4 from the corresponding powers in the right-hand-side of (2.13) . Thus, the proof follows as that of Theorem 2.1 in view of (2.18). □

It is obvious from Theorem 2.1 and Theorem 2.5 that for $m \in \mathbf{R}_+$ the sequences $\{1/x_k\}$ obtained from (2.1) and (2.13) are also positive and uniformly bounded satisfying $M^{-1} \le x_k^{-1} \le m^{-1}$. Formally, one gets the following results which follow directly by replacing $x_k \to 1/x_k$ in (2.1) and (2.13), respectively:

**Corollary 2.6**. Assume that Theorem 2.1 holds with $m \neq 0$. Then, the solution sequence of the discrete equation:

$$x_{k+1} = \frac{x_{k+1-\ell}^{p}}{A_k x_{k+1-\ell}^{p} + B_k x_{k+1-s}^{r}} \; ; \; \forall k \in \mathbf{Z}_{0+} \tag{2.19}$$

satisfies $0 < M^{-1} \le x_k \le m^{-1} < \infty$ provided that $M^{-1} \le x_k \le m^{-1}$; $\forall k \in \min(1-\ell, 1-s) \cap \mathbf{Z}$. □

**Corollary 2.7**. Assume that Theorem 2.5 holds with $m \neq 0$. Then, the solution sequence of the discrete equation

$$x_{k+1} = \frac{1}{A_k + \sum_{i=1}^{q} B_{ki} \dfrac{x_{k+1-s_i}^{r_i}}{x_{k+1-\ell_i}^{p_i}}} \; ; \; \forall k \in \mathbf{Z}_{0+} \tag{2.20}$$

satisfies $0 < M^{-1} \le x_k \le m^{-1} < \infty$ provided that $M^{-1} \le x_k \le m^{-1}$; $\forall k \in \min(1-\ell, 1-s) \cap \mathbf{Z}$. □

**Theorem 2.8**. The following properties hold:

**(i)** Assume that $p_i \in (0,1); \forall i \in \overline{q}$ and $a_m, a_M \in \mathbf{R}_+$. Then, all solution sequences of (2.13) generated under admissible bounded initial conditions are bounded.

**(ii)** Assume that $p \in (0,1)$ and $a_m, a_M \in \mathbf{R}_+$. Then, all solution sequences of (2.1) generated under admissible bounded initial conditions are bounded.

**Proof**: (i) The proof is made by contradiction. If a solution sequence $S = \{x_k\}$ of (2.13) is unbounded then there is a strictly monotone increasing subsequence $S_1 \subset S$ which diverges as $k \to \infty$ and which can be built as follows:

- $x_{k+1}$, $x_{k+1-\delta_k}$ are two consecutive elements of $S_1$ for $k \ge N_0 \in \mathbf{Z}_+$, $\delta_k \in \mathbf{Z}_+$ (being dependent on k ) for some $N_0$ finite and sufficiently large, such that for $\hat{p}_M := \max_{i \in \overline{q}} (p_i) \in (0,1)$

- $x_{k+1} = A_k + \sum_{i=1}^{q} B_{ik} \dfrac{x_{k+1-\ell_i}^{p_i}}{x_{k+1-s_i}^{r_i}} \le a_M + \sum_{i=1}^{q} B_{ik} \dfrac{x_{k+1-\ell_i}^{p_i}}{a_m^{r_i}} \le a_M + \sum_{i=1}^{q} B_{ik} \dfrac{x_{k+1-\delta_k}^{\hat{p}_M}}{a_m^{r_i}} \; ; \forall k \in \mathbf{Z}_{0+}$

Thus,

$$0 < v_k := x_{k+1}^2 - x_{k+1-\delta_k}^2 \le \overline{v}_k := a_M^2 - \left( 1 - \left( \sum_{i=1}^{q} \frac{B_{ik}}{a_m^{r_i} x_{k+1-\delta_k}^{1-\hat{p}_M}} \right)^2 - 2 \left( \sum_{i=1}^{q} \frac{B_{ik}}{a_m^{r_i - 1} x_{k+1-\delta_k}^{2-\hat{p}_M}} \right) \right) x_{k+1-\delta_k}^2$$

(2.21)

since $S_1$ is strictly monotone increasing. But, for some sufficiently large $\hat{x}(>1) \in \mathbf{R}_+$, there exists a sufficiently large finite $\mathbf{Z}_{0+} \ni N_0 = N_0(\hat{x})$ such that $x_{k+1-\delta_k} \ge \hat{x}$ and $\overline{v}_k < 0$ if $\mathbf{Z}_{0+} \ni k \ge N_0$ which is a contradiction to (2.21). Thus, all solutions of (2.13) are bounded for any admissible set of initial conditions and the proof of Property (i) is complete. Property (ii) follows in a similar way for the discrete sequence (2.1) which is a particular case of (2.13). □

*Remark 2.9.* A particular case of Theorem 2.8 (ii) for constant parametrical sequences $A_k = A$ and $B_k = 1$; $\forall k \in \mathbf{Z}_{0+}$, has been proven in [1]. □

### 3. Equilibrium points

The following two results follow from the limiting equations of (2.1), (2.13), (2.19) and (2.20) as the parameter sequences converge asymptotically to finite limits.

**Lemma 3.1**. Assume that $A_k \to A \in \mathbf{R}_+$, $B_{ik} \to B_i \in \mathbf{R}_{0+}$ as $k \to \infty$; $\forall i \in \overline{q}$ with $B_i \in \mathbf{R}_+$ for at least one $i \in \overline{q}$. Then, the following properties hold:

**(i)** $x \in \mathbf{R}$ is an equilibrium point of (2.13) if and only if

$$x - A = \sum_{i=1}^{q} B_i x^{r_i - p_i} \tag{3.1}$$

$x$ is an equilibrium point of (2.13) only if $x \ge A$.

$x = 1$ is an equilibrium point of (2.13) if and only if $A + \sum_{i=1}^{q} B_i = 1$.

**(ii)** $\overline{x} \in \mathbf{R}$ is an equilibrium point of (2.20) if and only if

$$\left( A + \sum_{i=1}^{q} B_i \overline{x}^{r_i - p_i} \right) \overline{x} = 1 \tag{3.2}$$

A necessary condition is $\overline{x} \le 1/A$. Also, $\overline{x} = 1$ is an equilibrium point of (2.20) if and only if $A + \sum_{i=1}^{q} B_i = 1$.

**(iii)** x= $\overline{x}$ = 1 is an equilibrium point of both (2.13) and (2.20) if and only if $A + \sum_{i=1}^{q} B_i = 1$.

**(iv)** Define the disjoint sets $Q_1:=\{i\in\overline{q}: r_i=p_i-1\}$, $Q_2:=\{i\in\overline{q}: r_i>p_i-1\}$, $Q_3:=\{i\in\overline{q}: r_i<p_i-1\}$. Then, $\overline{x}=0$ is the unique equilibrium point of (2.20) if and only if $r_i\ge p_i-1$; $\forall i\in\overline{q}$ (i.e. $Q_3=\varnothing$) and $\exists$(at least one) $i\in\overline{q}$ such that $p_i=r_i+1$ (i.e. $Q_1\ne\varnothing$) with $\sum_{i\in Q_1}B_i=1$.

**(v)** If $p_i=r_i+1$; $\forall i\in\overline{q}$ and $\sum_{i=1}^{q}B_i<1$ then $\overline{x}=\dfrac{1-\sum_{i=1}^{q}B_i}{A}$ is the unique equilibrium point of (2.20) so that $x=\dfrac{A}{1-\sum_{i=1}^{q}B_i}$ is then the unique equilibrium point of (2.13).

**(vi)** If

$$\sum_{i\in Q_1}B_i+\sum_{i\in Q_3}B_iA^{|r_i-p_i+1|}>1$$

then (2.20), and thus (2.13), have no equilibrium points.

**Proof**: Properties (i)-(ii) follows after direct substitution of $A_k\to A$, $B_{ik}\to B$, $x_k\to x$; $\forall k\in\mathbf{Z}_{0+}$ in (2.13) and (2.20), respectively. Note that $x_k\ge 0$ ($x_k>0$ if $A>0$; $\forall k\in\mathbf{Z}_{0+}$) from Theorem 2.5 what implies from (3.1)-(3.2) that any equilibrium point satisfies $x\ge A>0$ and $\overline{x}\le 1/A$, respectively, and note also that $x=\overline{x}=1$ in (3.1)-(3.2) if and only if $A+\sum_{i=1}^{q}B_i=1$ from Properties (i) –(ii) . Property (iii) has been proven. Note that $\overline{q}=Q_1\cup Q_2\cup Q_3$ so that (3.2) can be rewritten as:

$$A\overline{x}+\sum_{i\in Q_1}B_i+\sum_{i\in Q_2}B_i\overline{x}^{\,r_i-p_i+1}+\sum_{i\in Q_3}\frac{B_i}{\overline{x}^{|r_i-p_i+1|}}=1 \tag{3.3}$$

which is satisfied with $\overline{x}=0$ if and only if $Q_1\ne\varnothing$ with $\sum_{i\in Q_1}B_i=1$ and $Q_3=\varnothing$. The sufficiency part of Property (iv) has been proven. The necessary part is proven by contradiction as follows. Assume that $Q_3=\varnothing$ and $\sum_{i\in Q_1}B_i\ne 1$ with $Q_1\ne\varnothing$ or that $Q_1=\varnothing$ then (3.3) fails if $\overline{x}=0$. Then $Q_1\ne\varnothing$ with $\sum_{i\in Q_1}B_i=1$ is needed for $\overline{x}=0$ to be an equilibrium point of (2.20). Assume that $Q_3\ne\varnothing$. Then, the left –hand-side of (3.3) is $+\infty$ if $\overline{x}=0$ and thus (3.3) fails so that it is needed $Q_3=\varnothing$ for $\overline{x}=0$ to be an equilibrium point of (2.20). The uniqueness of the equilibrium point $\overline{x}=0$ follows from the uniqueness of the solution of the linear algebraic system $A\overline{x}+\sum_{i\in Q_1}B_i=1=\sum_{i\in Q_1}B_i$ obtained from (3.2). Property (v) has been proven. Note from (3.3) that an equilibrium point $x=\overline{x}^{-1}\ge A>0$ exists if and only if :

$$\sum_{i\in Q_2}\frac{B_i}{A^{r_i-p_i+1}}+\sum_{i\in Q_3}\frac{B_i}{\overline{x}^{|r_i-p_i+1|}}\ge 1-A\overline{x}-\sum_{i\in Q_1}B_i=\sum_{i\in Q_2}B_i\overline{x}^{\,r_i-p_i+1}+\sum_{i\in Q_3}\frac{B_i}{\overline{x}^{|r_i-p_i+1|}}$$

$$\geq \sum_{i\in Q_2} B_i \overline{x}^{\,r_i - p_i + 1} + \sum_{i\in Q_3} B_i A^{\left|r_i - p_i + 1\right|}$$

$$\Rightarrow 1 - \sum_{i\in Q_1} B_i - \sum_{i\in Q_3} B_i A^{\left|r_i - p_i + 1\right|} \geq \left( A + \sum_{i\in Q_2} B_i \overline{x}^{\,r_i - p_i} \right) \overline{x} \geq 0$$

$$\Rightarrow \sum_{i\in Q_1} B_i + \sum_{i\in Q_3} B_i A^{\left|r_i - p_i + 1\right|} \leq 1$$

$$\Leftrightarrow \left( \sum_{i\in Q_1} B_i + \sum_{i\in Q_3} B_i A^{\left|r_i - p_i + 1\right|} > 1 \Rightarrow \neg \exists\ \overline{x} \in \mathbf{R}_{0+} \text{ satisfying } (3.2) \right)$$

and Property (vi) has been proven. □

**Lemma 3.2**. Assume that $A_k \to A \in \mathbf{R}_+$, $B_k \to B \in \mathbf{R}_+$ as $k \to \infty$; $\forall i \in \overline{q}$ Then, the following properties hold:

**(i)** $x \in \mathbf{R}$ is an equilibrium point of (2.1) if and only if

$x - A = B x^{\,r-p}$

x is an equilibrium point of (2.1) only if $x \geq A$.

**(ii)** $\overline{x} \in \mathbf{R}$ is an equilibrium point of (2.19) if and only if

$\left( A + B\overline{x}^{\,r-p} \right) \overline{x} = 1$

A necessary condition is $\overline{x} \leq 1/A$. Also, $\overline{x} = 0$ is the unique equilibrium point of (2.19) if and only if $r = p - 1$ and B=1.

**(iii)** $x = \overline{x} = 1$ is an equilibrium point of both (2.1) and (2.19) if and only if $A + B = 1$.

**(iv)** If $p = r + 1$ and $B < 1$ then $\overline{x} = \dfrac{1-B}{A}$ is the unique equilibrium point of (2.20) so that $x = \dfrac{A}{1-B}$ is then the unique equilibrium point of (2.1).

If $p = r + 1$ and $B > 1$ then (2.19) has no equilibrium point so that (2.1) has no equilibrium point either.

**Proof**: The proofs of Properties (i)-(iii) are similar to that of Lemma 3.1 by the appropriate substitutions in (2.1) and (2.19), respectively. Property (iv) follows from Properties (iv) - (v) of Lemma 3.1, its last part being a particular case of Lemma 3.1 (vi). □

Some simple global stability results for any set of bounded initial conditions of the difference equations (2.1) and (2.13) follow in the sequel by excluding $\overline{x} := x^{-1} = 0$ to be an equilibrium point of the inverse sequences (2.19) –(2.20) , respectively. In the case that $\overline{x} = 0$ is a locally stable equilibrium point of (2.19) (respectively, (2.20)) then the solution of (2.1) ( respectively, (2.13) ) satisfies $x_k \to +\infty$ as $k \to \infty$ for certain sets of admissible bounded initial conditions even under the uniform boundedness results of Section 2 for different sets of bounded initial conditions.

**Proposition 3.3**. Assume that $Q_1 \neq \varnothing$, $Q_3 = \varnothing$ and $\sum_{i \in Q_1} B_i = 1$ .Then, there are unbounded solutions of (2.13) for bounded admissible initial conditions if the equilibrium point $\bar{x} = 0$ of (2.20) is locally stable.

**Proof**: $\bar{x} = 0$ is an equilibrium point of (2.20) from Lemma 3.1 (iv). Thus, if it is locally unstable, then there exist unbounded solutions of (2.13) for admissible bounded initial conditions. □

**Proposition 3.4**. If $\sum_{i \in Q_1} B_i + \sum_{i \in Q_3} B_i A^{|r_i - p_i + 1|} > 1$ then all the solutions of (2.13) and (2.20) are oscillatory and uniformly bounded, under admissible bounded initial conditions, in the sense that for each $k \in \mathbf{Z}_{0+}$, there exist finite $N_k, M_k \in \mathbf{Z}_+$, dependent in general on k, such that $\operatorname{sign}\left(x_{k+N_k+M_k} - x_{k+N_k}\right) = -\operatorname{sgn}\left(x_{k+N_k} - x_k\right)$.

**Proof**: From Lemma 3.1 (vi), the solution of (2.20), and then that of (2.13), have no equilibrium points. Thus, no solution of (2.20) converges asymptotically to zero so that no solution of (2.13) diverges asymptotically to $+\infty$. Thus, all solutions of (2.13) are bounded but non constant and not asymptotically constant (since there are no equilibrium points) for any admissible bounded initial conditions. As a result all the solutions are oscillatory uniformly bounded . □

**Proposition 3.5**. If $r = p - 1$ and B=1 then there are unbounded solutions of (2.19) for bounded admissible initial conditions if the equilibrium point $\bar{x} = 0$ of (2.1) is locally stable.

**Proof**: It is close to that of Proposition 3.3 from Lemma 3.2 (ii). □

**Proposition 3.6**. If $r = p - 1$ and $B > 1$ then all the solutions of (2.13) and (2.20) are oscillatory and uniformly bounded under admissible bounded initial conditions.

**Proof**: It is close to that of Proposition 3.4 from Lemma 3.2 (iv). □

It is clear that because of because of the summation of terms in the right- hand-side of (2.13) with a potential presence of possible rational powers, it is quite difficult to determine analytically the set of equilibrium points of that difference equation , provided it is non empty, and even the cardinal of this set . However, it is possible to determine their possible existence and upper and lower -bounds of their potential location by discussing the equilibrium solutions of upper- bounding and lower- bounding simpler difference equations provided that they exist. However, if the equilibrium points of such a lower-bounding and upper- bounding equations do not exist then it cannot be concluded the existence of equilibrium points of (2.13). To investigate those issues , first define:

$$\rho := \max_{i \in \bar{q}} \left(p_i - r_i\right) \ ; \quad \delta := \min_{i \in \bar{q}} \left(p_i - r_i\right) \tag{3.4}$$

so that $\delta \le \rho$. Assume that the real limits $A_k \to A$ ; $B_k \to B := \sum_{i=1}^{q} B_i$ as $k \to \infty$ exist and define the sequences:

$$\omega_{k+1} = A + B \min\left(\omega_k^{\delta}, \omega_k^{\rho}\right) \tag{3.5a}$$

$$v_{k+1} = A + B \max\left(v_k^{\rho}, v_k^{\delta}\right) \tag{3.5b}$$

$\forall k \in \mathbf{Z}_{0+}$ under initial conditions satisfying the constraints of (2.13). If the initial conditions of (3.5) are the same as any given admissible set of them for (2.13) then, for such a set:

$$\omega_k \le x_k \le v_k \,;\ \forall k \in \mathbf{Z} \tag{3.6}$$

The equilibrium points, if any, of the solutions of (3.5) are, respectively, the existing solutions of

$$f(\omega) := \omega - A \equiv g(\omega) := B \min\left(\omega_k^{\delta}, \omega_k^{\rho}\right) \,;\ f(v) \equiv h(v) := B \max\left(v_k^{\rho}, v_k^{\delta}\right) \tag{3.7}$$

*Remark 3.7*. In the particular case that $A \ge 1$ then $x_k \ge 1$; $\forall k(\ge k_o) \in \mathbf{Z}_{0+}$ for some finite $k_0 \in \mathbf{Z}_{0+}$ so that the equilibrium points, if any, of the solutions of (3.5) are, from (3.7), the respective existing solutions of

$$f(\omega) := \omega - A \equiv g(\omega) := B \omega_k^{\delta} \quad ; \ f(v) \equiv h(v) := B v_k^{\rho} \tag{3.8}$$

since (3.5) becomes:

$$\omega_{k+1} = A + B\,\omega_k^{\delta} \quad ; \quad v_{k+1} = A + B v_k^{\rho} \tag{3.9}$$

The result which follows is useful to allocate valid intervals for the allocation of the possible equilibrium points of (2.13), if any.

**Lemma 3.8**. Assume that $A \ge 1$. Then, the following properties hold:

**(i)** $f(\omega) \equiv g(\omega)$ has no solution if $\delta \ge 0$ so that (3.5a) has no equilibrium point.

**(ii)** $f(\omega) \equiv g(\omega)$ has a unique solution if $\delta < 0$ which satisfies the constraint:

$$A \le \omega = A + B\omega^{\delta} \le A + \frac{B}{A^{|\delta|}} \tag{3.10}$$

And which is the equilibrium point of (3.5a).

**(iii)** $f(v) \equiv h(v)$ has no solution if $\rho \ge 0$ (i.e. if $\exists i \in \overline{q} : p_i > r_i$) so that (3.5b) has no equilibrium point.

**(iv)** $f(v) \equiv h(v)$ has a unique solution if $\rho < 0$ (i.e. if $p_i < r_i ; \forall i \in \overline{q}$) which satisfies the constraint:

$$A \le v = A + B v^{\rho} \le A + \frac{B}{A^{|\rho|}} \tag{3.11}$$

Both difference equations (3.5) have respective unique equilibrium points which are solutions of $f(\omega) \equiv g(\omega)$ and $f(v) \equiv h(v)$, and satisfy the constraints :

$$1\le A\le \omega \le \min\left(v, A+B\omega^{\delta}\right)\le \min\left(v, A+\frac{B}{A^{|\delta|}}\right)\le 1+\frac{B}{A^{|\delta|}} \quad ; \quad 1\le A\le \omega\le v\le A+\frac{B}{A^{|\rho|}}\le 1+\frac{B}{A^{|\rho|}} \tag{3.12}$$

Both equilibrium points are identical if $\rho=\delta$.

**Proof**: **(i)** It follows from (3.8) since $f:\mathbf{R}_{0+}\to\mathbf{R}$ is linear with unity slope and $f(0)=-A$ and $g:\mathbf{R}_{0+}\to\mathbf{R}_{0+}$ is constant if $\delta=0$ and monotone strictly increasing growing faster than linearly if $\delta>0$ with $B=g(0)>f(0)=-A$. Thus, the functions $f:\mathbf{R}_{0+}\to\mathbf{R}$ and $g:\mathbf{R}_{0+}\to\mathbf{R}_{0+}$ have a unique coincidence point. **(ii)** It follows since $g:\mathbf{R}_{0+}\to\mathbf{R}_{0+}$ strictly monotone strictly decreasing growing faster than linearly if $\delta<0$, and $B=g(0)>f(0)=-A$ and $0=g(\infty)<f(\infty)=\infty$. Then, the functions $f:\mathbf{R}_{0+}\to\mathbf{R}$ and $g:\mathbf{R}_{0+}\to\mathbf{R}_{0+}$ have a unique coincidence point. Property **(iii)** and the first part of Property **(iv)** are proven in a similar way to Property **(i)** and Property **(ii)**, respectively. The second part of Property **(iv)** is a conclusion of Property **(ii)** and the first part of Property **(iv)** since $\omega\le v$, in view of (3.6), and $\rho=\max(\rho,\delta)<0\Rightarrow\delta<0$ provided that they exist but the last right-hand-side term of (3.10) is not necessarily upper-bounded by that of (3.11). □

The following result is direct from Lemma 3.1 and Lemma 3.8.

**Theorem 3.9**. Assume that $A\ge 1$ and that the empty or nonempty set of equilibrium points of (2.13) is $E=E_s\cup E_u$ where $E_s$ is the subset of locally asymptotically stable equilibrium points of (2.13) and $E_u$ is the subset of its locally unstable and locally critically stable equilibrium points.
Define nonempty sets of integers $Q_\delta:=\left\{i\in\overline{q}: r_i-s_i=\delta\right\}\subset\overline{q}$ and $Q_\rho:=\left\{i\in\overline{q}: r_i-s_i=\rho\right\}\subset\overline{q}$ ; and Define sets of $Q_\delta$ and $Q_\rho$ polynomials $T_{\delta j}(z)$ and $T_{\rho j}(z)$ of real coefficients as follows:

$$T_{\delta j}(z)=z^{n_{\delta j}}+\sum_{i=1}^{q}\left(\omega^{\delta-1}B_i\right)\left(p_j z^{n_{\delta j}-s_j}-r_j z^{n_{\delta j}-\ell_j}\right) \ , \ n_{\delta j}:=\max\left(s_j,\ell_j\right);\ \forall j\in Q_\delta$$

$$T_{\rho j}(z)=z^{n_{\rho j}}+\sum_{i=1}^{q}\left(v^{\rho-1}B_i\right)\left(p_j z^{n_{\rho j}-s_j}-r_j z^{n_{\rho j}-\ell_j}\right) \ , \ n_{\rho j}:=\max\left(s_j,\ell_j\right);\ \forall j\in Q_\rho$$

Then, the following properties hold:

**(i)** If $\rho<0$ then a unique $v$ exists so that $\overline{x}_i\in[A, v]$ ; $\forall\overline{x}_i\in E$ if $E\ne\varnothing$. The equilibrium point of (3.5b) is locally asymptotically stable if the $Q_\rho$ polynomials $T_{\rho j}(z)$ have all their zeros in the complex unit open disk centered at z=0, i.e. $C(0,1):=\left\{z\in\mathbf{C}:|z|<1\right\}$. Furthermore, $\limsup_{k\to\infty} x_k\le v$ for any set of bounded admissible initial conditions of (2.13) allocated within a sufficiently small ball centered at v.

**(ii)** If $s_i = \ell_i = 1; \forall i \in \bar{q}$, then the local asymptotic stability condition of the unique equilibrium point of Property (i) becomes in particular $v^{|\rho|+1} > |\rho|$ so that so that, in such a case, the equilibrium point of (3.5b) is asymptotically stable if it satisfies the subsequent constraint:

$$\left( A + \frac{B}{A^{|\rho|}} \right) \geq v > |\rho|^{1/(|\rho+1|)}$$

Furthermore, the solution sequence of (2.13) for bounded admissible initial conditions satisfies:

$$\max\left( A, |\rho|^{1/(|\rho+1|)} \right) \leq \limsup_{k \to \infty} x_k \leq \left( A + \frac{B}{A^{|\rho|}} \right)^{|\rho|+1}$$

If $|\rho|^{1/(|\rho+1|)} > A$ then the first of the above inequalities is strict.

**(iii)** If $\delta < 0$ then a unique $\omega$ exists so that $\bar{x}_i \in [\omega, \infty)$ if $E \neq \varnothing$. The equilibrium point of (3.5a) is locally asymptotically stable if the $Q_\rho$ polynomials $T_{\delta j}(z)$ have all their zeros in the complex unit open disk centered at z=0, i.e. $C(0,1) := \{ z \in \mathbf{C} : |z| < 1 \}$. Furthermore, $\liminf_{k \to \infty} x_k \geq \omega$ for any set of bounded admissible initial conditions of (2.13) allocated within a sufficiently small ball centered at $\omega$.

**(iv)** If $s_i = \ell_i = 1; \forall i \in \bar{q}$, then the local asymptotic stability condition of the unique equilibrium point of Property (iii) of the first difference equation in (3.5a) is locally asymptotically stable if $\omega^{|\delta|+1} > |\delta|$ so that, in such a case, it satisfies the subsequent constraint:

$$\left( A + \frac{B}{A^{|\delta|}} \right) \geq \omega > |\delta|^{1/(|\delta|+1)}$$

Furthermore, the solution sequence of (2.13) for bounded admissible initial conditions around $\omega$ satisfies:

$$\liminf_{k \to \infty} x_k \geq \max\left( \left( A + \frac{B}{A^{|\delta|}} \right)^{|\delta|+1}, |\delta|^{1/(|\delta|+1)} \right)$$

If $\left( A + \frac{B}{A^{|\delta|}} \right)^{|\delta|+1} < |\delta|^{1/(|\delta|+1)}$ then the above inequality is strict.

**(v)** $\max(\rho, \delta) < 0 \Rightarrow \bar{x}_i \in [\max(A, \omega), v]$ ; $\forall \bar{x}_i \in E$ if $E \neq \varnothing$ and the obtained bounds for $\limsup_{k \to \infty} x_k$ and $\liminf_{k \to \infty} x_k$ in Properties (i)-(iv) still hold.

**(vi)** $\left[ \left( p_i \leq r_{i+1} + 1; \forall i \in \bar{q} \right) \wedge \left( p_j = r_{j+1} + 1; \forall j \in Q_1 (\neq \varnothing) \subset \bar{q} \right) \right] \Rightarrow \left[ (E = \varnothing) \wedge \left( x_k \to \infty \quad \text{as} \quad k \to \infty \right) \right]$

provided that $\sum_{i \in Q_1} B_i = 1$.

**(vii)** $\left[\left(p_i = r_{i+1}+1; \forall j \in \bar{q}\right) \wedge \left(\sum_{i\in\bar{q}} B_i < 1\right)\right] \Rightarrow \left[\left(E \neq \emptyset\right) \supset \left\{ x = \frac{A}{1-\sum_{i=1}^{q} B_i} \right\}\right]$ .

**Proof**: The local asymptotic stability issues for the equilibrium points $\omega$ and v of the lower-bounding and upper-bounding solution sequences (3.5a) and (3.5b) for the solution of (2.13) follow directly from the linearized stability theorem [10] (see also [1] and [9]). Thus, the characteristic polynomials of the local perturbations $\Delta\omega$ and $\Delta v$ around such equilibrium points have to be stable assumed that those equilibrium points $\omega$ and v are locally asymptotically stable.

**Proofs of (i) and (iii)**: Consider the equilibrium point $v = \lim_{k\to\infty} v_k$ of (3.5b) and a local perturbation $\Delta v$ around it. Then, the local stability around such a point is characterized by the following dynamics:

$$\Delta v_{k+1} = \left(\sum_{i=1}^{q} B_i\right) \frac{p_j v^{r_j} v^{p_j-1} \Delta\omega_{k+1-\ell_i} - r_j v^{p_j} v^{r_j-1} \Delta v_{k+1-s_i}}{v^{2r_j}} \; ; \; \forall j \in Q_\rho$$

Define the one-step advance and delay operators $q$ and $q^{-1}$, respectively, as $x_{k+1} = q x_k$, $x_k = q^{-1} x_{k+1}$. Since $\rho < 0$ from Lemma 3.8 (iv) , the above first-order incremental equations around v , whose number is card $Q_\rho$, may be rewritten as $T_{\rho j}(z)\Delta v$, where $T_{\rho j}(z)$ is defined as follows for appropriately defined real coefficients $t_i$ from the polynomial identities below:

$$T_{\rho j}(z) = z^{n_{\rho j}} + \sum_{i=1}^{n_{\rho j}-1} t_{ij} z^{n_{\rho j}-i} \equiv z^{n_{\rho j}} + \sum_{i=1}^{q} \left(\frac{B_i}{v^{|\rho|+1}}\right)\left(p_j z^{n_{\rho j}-s_j} - r_j z^{n_{\rho j}-\ell_j}\right); \; \forall j \in Q_\rho$$

Then, v is locally asymptotically stable if and only if the polynomials $T_{\rho j}(z)$ ; $\forall j \in Q_\rho$ are all stable; i.e. with all their zeros lying in $|z| < 1$, where z is the complex indeterminate argument of the discrete z-transform (which is formally identical to the one-step-temporal advance operator q). Also, the boundedness relation (3.6) implies $\limsup_{k\to\infty} x_k \le v$ for initial conditions within a sufficiently small ball centered at v. The proof of Property **(iii)** for the equilibrium point $\omega$ of (3.5a) is similar to that of Property **(i)** with the replacements $\rho \to \delta$, $T_{\rho j}(z) \to T_{\delta j}(z)$, $n_{\rho j} \to n_{\delta j}$ , $Q_\rho \to Q_\delta$, $v \to \omega$ . Also, the boundedness relation (3.6) implies $\limsup_{k\to\infty} x_k \ge \omega$ for initial conditions within a sufficiently small ball centered at $\omega$.

**Proofs of (ii) and (iv)**: Now, $s_i = \ell_i = 1; \forall i \in \bar{q}$ and $\rho < 0$ or $\delta < 0$ . Then, $|\Delta v| < \left|\rho v^{\rho-1} \Delta v\right|$ irrespective of $\Delta v$ implying that $\omega^{|\rho|+1} > |\rho|$ so that v is locally asymptotically stable. Also, the limit superior of any solution of (2.13) for any set of bounded admissible conditions does not exceed the value of v.

$\left|\Delta\omega\right|<\left|\delta\omega^{\delta-1}\Delta\omega\right|$ irrespective of $\Delta\omega$ implying that $\omega^{|\delta|+1}>|\delta|$ so that $\omega$ is locally asymptotically stable. Also, the limit inferior of the solution of (2.13) for any set of bounded admissible initial conditions is non less than $\omega$. The remaining parts of the proof follow directly from Lemma 3.1 and Lemma 3.8.
The proof of **(v)** follows from the above results by noting that $\max(\rho,\delta)<0 \Leftrightarrow |\delta|\ge|\rho|$ since $\delta\le\rho$.
The proofs of (vi)-(vii) follow from Lemma 3.1 and Lemma 3.8. □

If $A<1$ then the above two lat results become modified as follows:

**Lemma 3.10**. Assume that $A<1$. Then, the following properties hold:
**(i)** $f(\omega)\equiv g(\omega)$ has no equilibrium point if $\delta\ge 0$.
**(ii)** $f(\omega)\equiv g(\omega)$ has a unique equilibrium point if $\delta<0$ which satisfies the constraint

$$A\le\omega=A+B\min\left(\omega^{\delta},\,\omega^{\rho}\right)\le A+B\min\left(\frac{1}{A^{|\delta|}},A^{\rho}\right)=A+BA^{\rho} \tag{3.13}$$

**(iii)** $f(v)\equiv h(v)$ has no equilibrium point if $\rho\ge 0$.
**(iv)** $f(v)\equiv h(v)$ has a unique equilibrium point if $\rho<0$ which satisfies the constraint

$$A\le v=A+B\max\left(v^{\delta},\,v^{\rho}\right)\le A+B\min\left(\frac{1}{A^{|\delta|}},\frac{1}{A^{|\rho|}}\right)\le A+\frac{B}{A^{|\delta|}} \tag{3.14}$$

**(v)** If $\max(\rho,\delta)<0$ then both difference equations (3.5) have unique equilibrium points which satisfy the constraints:

$$A\le\omega\le\min\left(v,A+\frac{B}{A^{|\rho|}}\right)\le A+\frac{B}{A^{|\delta|}}\le A+B;\quad A\le\omega\le v\le A+\frac{B}{A^{|\delta|}}\le A+B \tag{3.15}$$

Both equilibrium points are identical if $\rho=\delta$.

**Proof**: **(i)** It follows from (3.8) since $f:\mathbf{R}_{0+}\to\mathbf{R}$ is linear with unity slope and $f(0)=-A$ and $g:\mathbf{R}_{0+}\to\mathbf{R}_{0+}$ is constant if $\delta=0$ and monotone strictly increasing growing faster than linearly if $\delta>0$ with $B=g(0)>f(0)=-A$. Thus, the functions $f:\mathbf{R}_{0+}\to\mathbf{R}$ and $g:\mathbf{R}_{0+}\to\mathbf{R}_{0+}$ have a unique coincidence point. **(ii)** It follows since $g:\mathbf{R}_{0+}\to\mathbf{R}_{0+}$ strictly monotone strictly decreasing growing faster than linearly if $\delta<0$, and $B=g(0)>f(0)=-A$ and $0=g(\infty)<f(\infty)=\infty$. Then, the functions $f:\mathbf{R}_{0+}\to\mathbf{R}$ and $g:\mathbf{R}_{0+}\to\mathbf{R}_{0+}$ have a unique coincidence point. Property **(iii)** and Property **(iv)** are proven in a similar way to Property **(i)** and Property **(ii)**, respectively. Property **(v)** is a conclusion of Property **(ii)** and Property **(iv)** since $\omega\le v$, in view of (3.6), provided that they exist but the last right-hand-side term of (3.10) is not necessarily upper-bounded by that of (3.11). □

The following two results are direct from Lemma 3.1 and Lemma 3.8 for the case A<1.

**Theorem 3.11**. Assume that $A<1$ and that the empty or nonempty set of equilibrium points of (2.13) is $E=E_s \cup E_u$ where $E_s$ is the subset of locally asymptotically stable equilibrium points of (2.13) and $E_u$ is the subset of its locally unstable and locally critically stable equilibrium points. Then, Theorem 3.9 (i) and Theorem 3.9 (iii) still hold. □

**Theorem 3.12**. Assume that $A<1$. Then, if $s_i=\ell_i=1; \forall i \in \bar{q}$, the following properties hold:

**(i)** If $\rho<0$ then a unique $v$ exists so that $\bar{x}_i \in [A, v]$ ; $\forall \bar{x}_i \in E$ if $E \neq \varnothing$. The equilibrium point of the second difference equation in (3.5b) is locally asymptotically stable if $\bar{v}^{|\delta|+1}>|\delta|$ so that so that, in such a case, it satisfies the subsequent constraint:

$$\left(A+\frac{B}{A^{|\delta|}}\right) \geq v > |\delta|^{1/(|\delta|+1)}$$

**(ii)** If $\delta<0$ then a unique $\omega$ exists so that $\bar{x}_i \in [\omega, \infty)$ if $E \neq \varnothing$. The equilibrium point of the first difference equation in (3.5a) is locally asymptotically stable if $\omega^{|\rho|+1}>|\rho|$ so that, in such a case, it satisfies the subsequent constraint:

$$\left(A+\frac{B}{A^{|\rho|}}\right) \geq \omega > |\rho|^{1/(|\rho|+1)}$$

**(iii)** $\max(\rho, \delta)<0 \Rightarrow \bar{x}_i \in [\max(A, v), v]$ ; $\forall x_i \in E$ if $E \neq \varnothing$ and the obtained bounds for $\limsup_{k\to\infty} x_k$ and $\liminf_{k\to\infty} x_k$ in Properties (i)-(ii) still hold.

**(iv)** $[(p_i \leq r_{i+1}+1; \forall i \in \bar{q}) \wedge (p_j = r_{j+1}+1; \forall j \in Q_1(\neq \varnothing) \subset \bar{q})] \Rightarrow [(E=\varnothing) \wedge (x_k \to \infty \text{ as } k \to \infty)]$ provided that $\sum_{i\in Q_1} B_i = 1$.

**(v)** $$\left[(p_i = r_{i+1}+1; \forall j \in \bar{q}) \wedge \left(\sum_{i\in\bar{q}} B_i < 1\right)\right] \Rightarrow \left[(E \neq \varnothing) \supset \left\{ x = \frac{A}{1-\sum_{i=1}^{q} B_i} \right\}\right].$$

**Proof**: It is very close to that of Theorem 3.9 [(ii) , (iv), (v) , (vi) and ( vii) ] by noting that $\frac{1}{A^{|\delta|}} \geq \frac{1}{A^{|\rho|}}$ if $\delta \leq \rho < 0 \Leftrightarrow |\delta| \geq |\rho|$. Then, the proof is omitted. □

*Remark 3.13*. It is difficult to calculate the exact allocations equilibrium points of (2.13) for large values of q, except in simple cases, because the exponents in (2.13) are in general rational numbers. However, if such points can be calculated or approximated then the ideas in Theorem 3.9 and Theorem 3.11 about the use of the linearized stability theorem can be used to determine the local stability of each of those points. In particular, first-order incremental dynamic systems around the equilibrium points may be built for that purpose as follows:

$$\Delta x_{k+1} = \sum_{i=1}^{q} B_i \frac{p_i x^{r_i+p_i-1}\Delta x_{k+1-\ell_i} - r_i x^{r_i+p_i-1}\Delta x_{k+1-s_i}}{x^{2r_i}} \; ; \; \forall x \in E$$

This leads to the following polynomials that characterize the local stability properties of the incremental equations $T(x,z)\Delta x = 0 \, ; \forall x \in E$, where the polynomial $T(x,z)$ is parameterized at the equilibrium points of (2.13) and defined as follows:

$$T(x,z) = z^n + \sum_{i=1}^{n-1} t_i(x) z^{n-i} \equiv z^n + \sum_{i=1}^{q} B_i x^{p_i - r_i - 1}\left(p_i z^{n-s_i} - r_i z^{n-\ell_i}\right)$$

$$n := \max\left(\max_{i\in\overline{q}} \ell_i \, , \, \max_{i\in\overline{q}} s_i\right)$$

Thus, if $T(x,z)$ is stable then $x \in E$ is a locally asymptotically stable equilibrium point of (2.13). □

A practical test to investigate the local stability of the equilibrium points can be performed via comparison tests with known given stable polynomials as follows:

*Remark 3.14*. The stability character of the polynomials $T(x,z); \; \forall x \in E$ may be investigated from a sufficiency point of view by comparing them to given polynomials $T^*(x,z) := z^n + \sum_{i=1}^{n-1} t_i^*(x) z^{n-i}$ which are known to be stable as follows. Define error polynomials:

$$\tilde{T}(x,z) := T(x,z) - T^*(x,z) = \sum_{i=1}^{n-1} \tilde{t}_i(x) z^{n-i}$$

$$= \sum_{i=1}^{n-1}\left(t_i(x) - t_i^*(x)\right) z^{n-i} = \sum_{i=1}^{q} B_i x^{p_i - r_i - 1}\left(p_i z^{n-s_i} - r_i z^{n-\ell_i}\right) - \sum_{i=1}^{n-1} t_i^*(x) z^{n-i}$$

$; \; \forall x \in E$. Define indicator subsets of $\overline{q}$:

$$Q_{si} := \left\{ i \in \overline{q} : n - s_i = i \right\} \; ; \quad Q_{\ell i} := \left\{ i \in \overline{q} : n - \ell_i = i \right\} \; ; \; \forall i \in \overline{q}$$

which can be empty. Thus, the coefficients of the error polynomials are given by:

$$\tilde{t}_i(x) = \left( \sum_{j \in Q_{si}} p_j B_j - \sum_{j \in Q_{\ell i}} r_j B_j \right) x^{p_i - r_i - 1} - t_i^*(x)$$

The Rouché theorem for zeros of analytic functions establishes that $T(x,z)$ is stable (i.e. with all its zeros in the complex open unit disk $|z|<1$) if $T^*(x,z)$ is stable and, furthermore, $\left|\tilde{T}(x,z)\right| < \left|T^*(x,z)\right|$ on the boundary of the unit disk $|z|=1$. In fact, such a property guarantees that $T(x,z)$ has the same number of zeros in $|z|<1$ than $T^*(x,z)$ so that $T(x,z)$ is stable since $T^*(x,z)$ is also stable for any given $x \in E$ (see, for instance, [11], [12]). □

**ACKNOWLEDGMENTS**

The author is grateful to the Spanish Ministry of Education by its partial support of this work through Grant DPI2006-00714. He is also grateful to the Basque Government by its support through Grants GIC07143-IT-269-07and SAIOTEK S-PE08UN15.

**REFERENCES**

[1] X. Yang, F. Sun and C. Zhang, "On the recursive sequence $x_n = A + \frac{x_{n-k}^p}{x_{n-1}^r}$", *Discrete Dynamics in Nature and Society*, Vol. 2009 , Article ID in press.

[2] A.M. Amleh, E.A. Grove, and G. Ladas, "On the recursive sequence $x_{n+1} = \alpha + \frac{x_{n-1}}{x_n}$", *Journal of Mathematical Analysis and Applications*, Vol. 233, pp. 790-798, 1999.

[3] K.S. Berenhaut and S. Stevic, "On positive nonoscillatory solutions of the difference equation $x_{n+1} = \alpha + \frac{x_{n-k}^p}{x_n^p}$", *Journal of Difference Equations and Applications*, Vol. 12, No. 5, pp.495-499, 2006.

[4] K.S. Berenhaut and S. Stevic, "The behavior of the positive solutions of the difference equation $x_n = A + \left(\frac{x_{n-2}}{x_{n-1}}\right)^p$", *Journal of Difference Equations and Applications*, Vol. 12, No. 9, pp.909-918, 2006.

[5] S. Stevic, "On the recursive sequence $x_{n+1} = \alpha + \frac{x_{n-1}^p}{x_n^p}$", *Journal of Applied Mathematics and Computing*, Vol. 18, Nos. 1-2, pp.229-234, 2005.

[6] S. Stevic, "On the recursive sequence $x_{n+1} = \alpha + \frac{x_n^p}{x_{n-1}^p}$", *Discrete Dynamics in Nature and Society*, Vol. 2007, Article ID 34517, 9 pages.

[7] S. Stevic, "On the recursive sequence $x_{n+1} = \alpha + \frac{x_n^p}{x_{n-1}^r}$", *Discrete Dynamics in Nature and Society*, Vol. 2007, Article ID 40963, 9 pages.

[8] M. Dehghan and R. Mazrooei-Sebdani, "Some characteristics of solutions of a class of rational difference equations", *Kybernetes*, Vol. 37, No. 6, pp. 786-796, 2008.

[9] J.S. Yu, " Asymptotic stability for a linear difference equation with variable delay", *Computers and Mathematics with Applications*, Vol. 36, Nos. 10-12, pp. 203-210, 1998.

[10] M.R.S. Kulenovic and G. Ladas, *Dynamics of a Nonlinear Difference Equation with Open Problems and Conjectures*, Chapman and Hall /CRC, Boca Raton , Fl. , 2002.

[11] M. De la Sen , " Sufficiency- type stability and stabilization criteria for linear time-invariant systems with constant point delays", *Acta Applicandae Mathematicae*, Vol. 83, No. 3, pp. 235-256, 2004.

[12] M. De la Sen, " Stabilization criteria for continuous-time linear time-invariant systems with constant lags", *Discrete Dynamics in Nature and Society*, Vol. 2006, Article ID 87062, 19 pages, doi:10.1155/DDNS/2006/87062.

[13] G. Papaschinopoulos, G. Stefanidou and C.J. Schinas, " Boundedness, attactivity, and stability of a rational difference equation with two periodic coefficients", *Discrete Dynamics in Nature and Society*, Vol. 2009, Article ID 973714, 19 pages, doi:10.1155/DDNS/2009/973714.

[14] I. Yalcinkaya, "On the difference equation $x_{n+1} = \alpha + x_{n-m} / x_n^k$", *Discrete Dynamics in Nature and Society*, Vol. 2008, Article ID 805460, 8 pages, doi:10.1155/DDNS/2008/805460.

[15] I. Yalcinkaya, " On the global asymptotic stability of a second-order system of difference equations", *Discrete Dynamics in Nature and Society*, Vol. 2008, Article ID 860152, 12 pages, doi:10.1155/DDNS/2008/860152.

[16] M. De la Sen and NS Luo, "On the uniform exponential stability of a wide class of linear time-delay systems", *Journal of Mathematical Analysis and Applications*, Vol. 289, No. 2, pp. 456-476, 2004.